\documentclass[12pt]{amsart}
\usepackage{verbatim}
\usepackage{eucal}
\usepackage{amssymb}
\usepackage{mathrsfs}
\usepackage{graphicx}
\usepackage{psfrag}
\newif \ifwide
\newif \ifavnermargin
\def \makemargins{
\ifwide
        \oddsidemargin .25in
        \evensidemargin .25in
        \textwidth 6.00in
\else
\fi
\ifavnermargin
        \headheight=7pt
        \textheight=574pt
        \textwidth=432pt
        \topmargin=14pt
        \oddsidemargin=18pt
        \evensidemargin=18pt
\else   
\fi
}

\avnermargintrue
\makemargins

\theoremstyle{plain}
\newtheorem{theorem}[subsection]{Theorem}
\newtheorem{proposition}[subsection]{Proposition}
\newtheorem{corollary}[subsection]{Corollary}

\theoremstyle{definition}
\newtheorem{definition}[subsection]{Definition}
\newtheorem{example}[subsection]{Example}

\theoremstyle{remark}
\newtheorem{question}[subsection]{Question}

\newcommand{\PP}{{\mathbb P}}

\newcommand{\CC}{{\mathbb C}}

\renewcommand{\theta}{\vartheta}

\begin{document}

\title{Wonderful blowups associated to group actions}

\newif \ifsections

\sectionstrue

\def \makeauthor{
\author{Lev A. Borisov}
\address{Department of Mathematics\\
Columbia University\\
New York, NY  10027}
\email{lborisov@math.columbia.edu}

\author{Paul E. Gunnells}
\address{Department of Mathematics and Computer Science\\
Rutgers University\\
Newark, NJ  07102--1811
}
\email{gunnells@andromeda.rutgers.edu}
}

\makeauthor

\date{September 1, 2000}

\thanks{The second author was partially supported by the NSF}

%
%

\begin{abstract}
A group action on a smooth variety provides it with the natural
stratification by irreducible components of the fixed point
sets of arbitrary subgroups. We show that the corresponding maximal 
wonderful blowup in the sense of MacPherson-Procesi has only
abelian stabilizers. The result is inspired by the abelianization
algorithm of Batyrev.
\end{abstract}
\maketitle

%
%
\section{Introduction}\label{introduction}

\subsection{}
Let $X/\CC$ be a smooth variety, and for $n>1$ consider the product
$X^{n}$. For any collection $I=\{I_k\}$ of disjoint subsets
$I_k\subseteq \{1,\dots ,n \}$ with each $I_k$ of order at least two, let
$\Delta_{I}$ be the polydiagonal consisting of all $x\in X^{n}$ with
$x_{i}=x_{j} $ if and only if $i,j\in I_k$ for some $k$.  The locally
closed subvarieties $\Delta_I$ provide a stratification of $X^{n}$
that is stable under the permutation action of the symmetric group
$S_{n}$.  
In \cite{fm}, Fulton and MacPherson constructed a
smooth variety $X[n]$ and a birational morphism $\pi \colon
X[n]\rightarrow X^{n}$ such that $\pi $ is an isomorphism over
the open set $\Delta_{\varnothing}$, and such that the
stabilizers of the induced $S_{n}$ action on $X[n]$ are solvable. More
recently, Ulyanov \cite{Ulyanov} constructed a space $X\langle n
\rangle$ with a stratified birational morphism to $X[n]$ such that the
stablilizers of the induced $S_{n}$ action are abelian.  The variety
$X\langle n \rangle$ is constructed by systematically blowing up the
closures of the preimages of polydiagonals $\Delta_{I}$ of increasing
dimension.

Both $X[n]$ and $X\langle n \rangle$ are special cases of a general
construction of MacPherson and Procesi \cite{MacPh.Pr}.  Suppose that
$Z$ is an irreducible smooth variety and that $S=\{S_i\}$ 
is a set of locally closed  subvarieties endowing
$Z$ with a stratification.  Suppose further that this stratification
is \emph{conical}; this means that every point $z\in Z$ has an
analytic neighborhood isomorphic to a product of a trivial stratification
of a disc and of a restriction to a disc of a $\CC^*$-invariant
stratification of a vector space.  A basic example is $Z=X^{n}$ and
$S=\{\Delta_I\}$.  Then MacPherson and Procesi constructed a family of
\emph{wonderful} compactifications of the open stratum.  These
smooth compactifications map birationally onto $Z$ isomorphically away
from strata of codimension bigger than one, and have the property that
the complement of the open stratum is a normal crossing divisor.
Moreover, for any $S$ there
are canonical minimal and maximal wonderful compactifications.  In the
example above, these spaces are respectively $X[n]$ and $X\langle n
\rangle$.

\subsection{}
In this paper we explain the connection
\begin{align*}
X[n]&\Longleftrightarrow \hbox{$S_{n}$ acts with solvable stabilizers}\\
X\langle n\rangle&\Longleftrightarrow \hbox{$S_{n}$ acts 
with abelian stabilizers}
\end{align*}
by investigating wonderful blowups associated to certain stabilizer
stratifications of a variety with group action.  More precisely, let
$X$ be a smooth variety and let $G$ be a finite group acting on $X$.
We define two stratifications of $X$, the \emph{stabilizer} and $Y$
stratifications, and show that they are conical.  Our main results
(Corollary \ref{abelianY}, Theorem \ref{main}, and Corollary
\ref{minYcor}) are that $G$ acts on
the maximal wonderful blowups associated to each of these
stratifications with abelian stabilizers, and on the minimal wonderful
blowup associated to the $Y$ stratification with solvable stabilizers.
The main tool is the abelianization algorithm of Batyrev
(\cite{bat.paper}, corrected in \cite{bat.correction}).

\subsection{Notations}
We fix a finite group $G$ for the rest of this paper, and work in the
category of smooth $G$-varieties over $\CC$, or smooth complex
$G$-manifolds.  We always assume that the action of $G$ is effective,
i.e. that only the identity stabilizes all of any $G$-variety $X$.  If
$X$ is a $G$-variety, then for every point $x\in X$ its stabilizer is
denoted by ${\rm Stab}(x)$. For every subgroup $H\subseteq G$ the set
of fixed points of $H$ is denoted by ${\rm Fixed}(H)$.  The latter is
always a disjoint union of smooth subvarieties of $X$.

\section{Batyrev's abelianization algorithm}\label{s1}
\subsection{}
The goal of this section is to describe a version of Batyrev's
abelianization algorithm 
\cite{bat.paper,bat.correction}. 
Our version is a bit easier
to formulate but generally results in more blowups.

Let $X$ be a $G$-variety as above. For every point $x\in X$ and any subgroup
$H$ fixing $x$ ($H$ may be a proper subgroup of ${\rm Stab}(x)$)
consider the decomposition of the tangent space $TX_x$ into irreducible
$H$-modules
$$
TX_x = \bigoplus_{i=1}^k V_i,
$$
where $\dim(V_i)=1$ for $i\leq k_1$ and $\dim(V_i)>1$ for
$i>k_1$. Each one-dimensional representation $V_i$ gives a character
$H\rightarrow \CC^*$. Denote by $H_1$ the common kernel of these
characters, and let $Y(x,H)$ be the subvariety of $X$ that is the
connected component of ${\rm Fixed}(H_1)$ passing through $x$. Observe
that the tangent space to $Y(x,H)$ at $x$ is precisely
$\oplus_{i=1}^{k_1}V_i$ \cite[Lemma 3]{bat.correction}. Notice that
$Y(x,H)=X$ if and only if $H$ is abelian. Otherwise, the codimension
of $Y(x,H)$ is at least $2$.  We will consider the set of the {\em
proper} subvarieties $Y_j\subseteq X$ that are equal to $Y(x,H)$ for
at least one pair $(x,H)$.

\begin{theorem}\label{batyrev.algorithm}\cite{bat.paper,bat.correction}
In the above setup, denote by $r(X)$ the maximum number of different
subvarieties $Y_j$ with non-trivial intersection.  Denote by $Z$ the
set of points $x\in X$ that are contained in $r(X)$ different
$Y_j$. Then $Z$ is smooth and $G$-equivariant.  We let $X_1$ be the
blowup of $X$ along $Z$, and iterate the procedure: for $i\geq 1$ we
compute $r (X_{i})$ and $Z_{i}$, and define $X_{i+1}$ to be the blowup
of $X_{i}$ along $Z_{i}$.  We claim that this process terminates after
at most $r(X)$ steps, and that the stabilizer of every point of the
resulting variety is abelian.
\end{theorem}

\begin{proof}
It is straightforward to see that the intersection of any number 
of the $Y_j$ is smooth. This implies that $Z$ is smooth, and it is 
clearly $G$-equivariant. It will now suffice to prove that $r(X_1)<r(X)$.

We claim that every $Y_{j,1}$ for $X_1$ is a proper preimage of 
some $Y_j$ for $X$ that is not contained in $Z$ and vice versa.
Let $Y_{j,1}=Y(x_1,H)$ be a proper subvariety of $X_1$. Denote 
the image of $x_1$ under the blowdown map by $x$, and consider
$Y(x,H)$. If $x\notin Z$, there is nothing to prove. If $x\in Z$,
notice that $Z$ is a subvariety of $Y=Y(x,H)$, so the tangent 
space $TZ_x$ is an $H$-submodule of
$$
TY_x=\bigoplus_{i=1}^{k_1}V_i.
$$
If the connected component of $Z$ passing through
$x$ equals $Y$, then the fiber of the blowup over $x$, which 
equals $\PP(\oplus_{i>k_1}V_i)$, has no $H$-fixed points.
Therefore the tangent space to $Z$ is a proper submodule of 
$TY_x$, which we assume to equal $\oplus_{i=1}^{k_2}V_i$.
The point $x_1$ then corresponds to a one-dimensional
submodule of $TY_x/TZ_x$, which we will identify with $V_{k_2+1}$.
The tangent space $(TX_1)_{x_1}$ is isomorphic as an $H$-module 
to
$$
\bigoplus_{i=1}^{k_2+1}V_i \bigoplus_{i=k_2+2}^{k}(V_{k_2+1}\otimes V_i).
$$
Therefore, the tangent space to $Y_{j,1}$ at $x_1$ is
$$
\bigoplus_{i=1}^{k_2+1}V_i \bigoplus_{i=k_2+2}^{k_1}(V_{k_2+1}\otimes V_i),
$$
which is the tangent space of the proper preimage of $Y$.
The statement in the opposite direction is proved similarly.

To show that $r(X_1)<r(X)$, suppose that $x_1$ is contained
in $r(X)$ different proper preimages of the $Y_j$.  Then $x_{1}$ corresponds
to a normal direction to $Z$ which is contained in all these $Y_j$'s,
which leads to a contradiction.
\end{proof}

\section{Batyrev's abelianization as a maximal wonderful blowup}
\subsection{}
The proof of Theorem \ref{batyrev.algorithm} shows that once the $Y_j$ 
are defined for $X$, their exact nature is not important for 
the further blowups. Indeed, at every step one blows up the intersection
set of the maximum number of the birational proper preimages of the $Y_j$.
This allows us to identify the result of Batyrev's algorithm
with the \emph{maximal wonderful blowup} of the MacPherson-Procesi family 
of blowups \cite{MacPh.Pr} associated to the stratification defined 
by the $Y_j$. 

\begin{definition}
The \emph{$Y$ stratification} of $X$ is the stratification induced by
the $Y_j$ as follows. For every subset of $\{Y_j\}$ one defines a
stratum that consists of all points that lie in all $Y_j$ from
the subset, but do not lie in any other $Y_j$. The empty strata
are then ignored.
\end{definition}

\begin{corollary}\label{abelianY}
The $Y$ stratification of $X$ is conical.  The maximal wonderful
blowup coincides with the result of Batyrev's algorithm above.
In particular, it has only abelian stabilizers.
\end{corollary}

\begin{proof}
First, it is easy to see that the $Y$ stratification is conical.
Moreover, this can be said about any stratification induced by some
connected components of fixed point sets of some subgroups of $G$.
Indeed, for any $x\in X$ there exists a ${\rm Stab}(x)$-equivariant 
isomorphism of a neighborhood $U\ni x$  
and a neighborhood of the origin in $TX_x$. Under this isomorphism
all strata map to linear subspaces, so the stratification is 
conical.

Second, the maximal wonderful blowup is defined by successively
blowing up proper preimages of all strata, starting with the strata of
smallest dimension.  This is equivalent to blowing up strata in any
order that is compatible with the partial ordering on the strata, and
Batyrev's algorithm clearly provides that.
\end{proof}

\ifsections
\section{The stabilizer stratification}
\subsection{}
\else
\subsection{The stabilizer stratification}
\fi

The $Y$ stratification is not local, in the sense
that its construction does not commute with $G$-equivariant 
open embeddings.
A more natural stratification of $X$ is induced by the 
set of all connected components of ${\rm Fixed}(H)$ for all subgroups
$H\subseteq G$. We will call it {\em stabilizer stratification} of $X$.
It is also conical by the argument presented in the proof of Corollary
\ref{abelianY}. 

\begin{proposition}\label{maxtoab}
Let $X^{Y}$ be the maximal wonderful blowup associated to the
$Y$ stratification of $X$, and let $X^{stab}$ be the maximal wonderful
blowup associated to the stabilizer stratification of $X$.  Then there
exists a natural $G$-equivariant map
$$
X^{stab}\longrightarrow  X^{Y}.
$$
\end{proposition}

\begin{proof}
The statement is clearly local in $X$, so it could be assumed that 
$X$ is a neighborhood of the origin in a $G$-vector space $V$. 
Denote by $S_Y$ and $S_{stab}$ the sets of strata in the $Y$ and
stabilizer stratifications respectively.
According to \cite[Definition p.461, Proposition p.470]{deConcini.Procesi} 
the maximal blowups 
$X^Y$ and $X^{stab}$ can be described as closures of the images of an
open subset of $V$ under the map to 
$$V\times \prod_{k\in S}\PP(V/V_k).
$$
Here $V_k$ is the tangent space to the stratum, and $S$ is 
either $S_Y$ or $S_{stab}$. The functoriality is then obvious.
\end{proof}

As a corollary we get the following theorem, which is the main result
of this paper.

\begin{theorem}\label{main}
The stabilizer of every point of the maximal wonderful blowup 
associated to the stabilizer stratification is abelian.
\end{theorem}

\begin{proof} Combine Corollary
\ref{abelianY} and Proposition \ref{maxtoab}.
\end{proof}

\ifsections
\section{Solvable stabilizers}\label{solvable}
\subsection{}
\else
\subsection{Solvable stabilizers}\label{solvable}
\fi
It is generally much easier to ensure that a blowup of a $G$-variety
has only {\em solvable} stabilizers, due to the following observation.

\begin{proposition}\label{criterion}
Let $\pi\colon X_1\rightarrow X$ be a $G$-equivariant birational
morphism of $G$-varieties.  Then all stabilizers of $X_1$ are solvable
if and only if the image of the exceptional divisors of $\pi$ contains
all points of $X$ with non-solvable stabilizers.
\end{proposition}

\begin{proof}
Suppose $X_1$ has points with non-solvable stabilizers. Consider the
point $x_1\in X_1$ with a non-solvable stabilizer $H$ that lies in the
minimum number of exceptional divisors of $\pi$. If this minimum
number is zero, then $X$ has a point with non-solvable stabilizer that
lies outside the image of the exceptional divisors.  Otherwise, let
$E$ be an exceptional divisor of $\pi$ that passes through $x$. The
tangent space $TX_{x_1}$ splits as an $H$-module into $TE_{x_1}$ and a
one-dimensional module. This one-dimensional module corresponds to a
character $\chi\colon H\rightarrow \CC^*$. The kernel of this
character is again non-solvable; however its fixed point set contains
points in the neighborhood of $x_1$ not lying in $E$. This
contradicts the minimality of $x_1$ and proves the ``if'' part. The
``only if'' part is obvious.
\end{proof}

\begin{corollary}\label{minYcor} 
The minimal wonderful blowup associated to the $Y$ stratification
has only solvable stabilizers.
\end{corollary}

\begin{proof}
In view of Proposition \ref{criterion}, one needs to show that every
point $x\in X$ with non-solvable stabilizer is contained in an
irreducible stratum of codimension bigger than one. Since all $Y_j$
have codimension at least two, it remains to observe that each point
$x$ with a non-solvable stabilizers is contained in $Y=Y(x,{\rm
Stab}(x))\not =X$.
\end{proof}

\begin{question}
Is it true that all stabilizers of the minimal wonderful blowup associated
with the stabilizer stratification are solvable?
\end{question}

This is really a question of whether there exists a non-solvable group
$G$ acting on $\CC^d$ with a fixed basis, such that for every $g\in G$
the space ${\rm Ker}(g-1)$ is a coordinate subspace. 

\section{Examples}
\begin{example}
As in Section \ref{introduction}, let $G$ be the symmetric group $S_n$
acting on the product $X^{n}$ of $n$ copies of a variety $X$. Then the
$Y$ stratification differs from the stabilizer stratification only in
that it does not distinguish general points on the large diagonals
(the conjugates of the preimages of the diagonals in $X\times X$) from the
general points on $X^n$. As the result, the maximal wonderful blowup
associated to the stabilizer stratification \cite{Ulyanov} is
obtained from the maximal wonderful blowup for the $Y$ stratification
by blowing up disjoint proper preimages of the large diagonals. 
Theorem \ref{main} implies the result of Ulyanov \cite{Ulyanov} that 
the stabilizer of the maximal blowup is abelian.
It is also easy to see that the discussion of Section
\ref{solvable} implies that the stabilizers of the minimal wonderful
blowups are solvable, which was first observed in \cite{fm}.
\end{example}

It is natural to ask whether for any finite group action of $G$ on $X$
there exists a \emph{minimum abelianization}, defined to be the
$G$-equivariant birational morphism $\pi\colon X_1\rightarrow X$ such
that $X_1$ has abelian stabilizers and such that every other such map
factors through $\pi$. Unfortunately, the next example shows that this
is not the case.

\begin{example}
Let $X=\CC^4$ and let $G=S_3$ act on $X$ by linear transformations
such that $X$ is a sum of two irreducible $2$-dimensional modules of
$G$.  Then a blowup at the origin or a blowup along any
$2$-dimensional submodule of $X$ has only abelian stabilizers, but
there is clearly no minimum abelianization.
\end{example}

%
%

\bibliographystyle{amsplain}

\providecommand{\bysame}{\leavevmode\hbox to3em{\hrulefill}\thinspace}

\end{document}